%% file: 2001-12.tex
\documentclass{gtart}

\input gtoutput

\lognumber{160}
\volumenumber{5}\papernumber{12}\volumeyear{2001}
\pagenumbers{347}{367}
\accepted{6 April 2001}
\received{09 January 2001}
\revised{10 March 2001}
\published{6 April 2001}
\proposed{Walter Neumann}
\seconded{David Gabai, Joan Birman}

\usepackage{rlepsf}
\let\relabela\adjustrelabel

\def\a{\alpha}
\def\b{\beta}
\def\g{\gamma}
\def\d{\delta}

\def\s{\sigma}

\def\Hull{{\rm Hull}}

\newtheorem{theorem}{Theorem}[section]
\newtheorem{lemma}[theorem]{Lemma}
\newtheorem{corollary}[theorem]{Corollary}

\newtheorem{proposition}[theorem]{Proposition}

\theoremstyle{remark}

\begin{document}
\title{Some surface subgroups survive surgery}
\authors{D Cooper and D\thinspace D Long}
\asciiauthors{D Cooper and D D Long}
\coverauthors{D Cooper and D\noexpand\thinspace D Long}

\address{Department of Mathematics, University of California\\Santa
Barbara, CA 93106, USA}
\email{cooper@math.ucsb.edu\\long@math.ucsb.edu}

\begin{abstract}
It is shown that with finitely many exceptions, the fundamental group
obtained by Dehn surgery on a one cusped hyperbolic 3--manifold contains
the fundamental group of a closed surface.
\end{abstract}

\asciiabstract{
It is shown that with finitely many exceptions, the fundamental group
obtained by Dehn surgery on a one cusped hyperbolic 3-manifold contains
the fundamental group of a closed surface.
}

\keywords{Dehn surgery, surface subgroups}

\primaryclass{57M27}

\secondaryclass{57M50, 20H10}

\maketitle

\section{Introduction}
A central unresolved question in the theory of closed hyperbolic
$3$--manifolds is whether they are covered by manifolds which contain
closed embedded incompressible surfaces.  An affirmative resolution of this
conjecture would imply in particular that all closed hyperbolic
$3$--manifolds contain the fundamental group of a closed surface of genus
at least two.  Even the simplest case of this conjecture, namely that of
the manifolds obtained by surgery on a hyperbolic manifold with a single
cusp has remained open for many years.  In this article we prove the
following theorem:
\begin{theorem}
\label{surgerytheorem}
Suppose that $M$ is a hyperbolic $3$--manifold with a single torus
boundary component.

Then all but finitely many surgeries on $M$ contain the fundamental group
of a closed orientable surface of genus at least two.
\end{theorem}
Our proof rests upon:
\begin{theorem}
\label{main}
Suppose that $S$ is an incompressible, $\partial$--incompressible
quasifuchsian surface with boundary slope $\a$.

Then there is a $K > 0 $ so that if $\g$ is any simple curve on $\partial M$
with $\Delta(\a,\g) > K$, the Dehn filled manifold $M(\g)$ contains the
fundamental group of a closed surface of genus at least two.
\end{theorem}
This result is similar in spirit to the main theorem of \cite{CL}; that
result applied only to surfaces of slope zero, however in that context we
were able to give an explicit (and fairly small) value for
$K$.

\proof[Proof of \ref{surgerytheorem} from \ref{main}]
It follows from \cite{CG} that $M$ contains at least two distinct strict boundary
slopes; and we shall show using Proposition 1.2.7 loc.~cit.\  (see
Lemma \ref{qf_representative}) that both these slopes are represented by
quasifuchsian surfaces.  Then \ref{main} implies the result. 
\endproof

 We now outline the proof.  Using the quasifuchsian
surface, we are able to construct a certain complex in the universal
covering of $M$.  This in turn gives rise to a map of a manifold with
convex boundary $f \co  X \rightarrow M$.  We are able to prove that this map
extends to a map of Dehn filled manifolds, $f_\g \co  X(\g) \rightarrow
M(\g)$.  The explicit nature of $X(\g)$ makes it possible to prove that
$X(\g)$ contains a surface group: a geometrical argument using convexity
then shows that $f_{\g*}$ is $\pi_1$--injective, completing the proof. 
The intuition for where the surface comes from in this construction is the
same as \cite{CL}: two surfaces glued together by a very long annulus
should remain incompressible in fillings distant from the boundary slope. 
It seems worth pointing out that it follows that the surface group which we
produce comes from an immersion of a surface into $M(\g)$ without triple
points.

A proof of Theorem \ref{main} by completely different methods has been
given by T Li (see \cite{Li}). His proof has the advantage that it
gives bounds on the constant $K$.

\section{The proof of \ref{main}}
Throughout this article we fix a hyperbolic manifold $M$ with a single
torus boundary component.  This gives rise to a discrete, faithful
representation $\rho\, \co  \pi_1(M) \longrightarrow PSL(2,{\bf C})$ which is
unique up to conjugacy.
\subsection{Quasifuchsian surfaces}

We recall that if $f\co  S \longrightarrow M$ is a proper map which is
injective at the level of fundamental groups then the surface $f(S)$ is
said to be {\em quasifuchsian } if the limit set of the group $\rho
f_*(\pi_1(S)))$ is a topological circle.  Many equivalent formulations
exist; for our purpose, it suffices to observe that it is shown in Chapter
V, Corollary 9.2 of \cite{BM} that a representation of a surface group is
quasifuchsian if and only if the representation is geometrically finite and
contains no accidental parabolics.  (We recall that a surface is said to
have no {\em accidental parabolics} if the conjugacy classes of elements
representing boundary components in $\pi_1(S)$ are exactly the conjugacy
classes that are parabolic under the representation.)

We begin with a purely topological lemma.
\begin{lemma}
\label{no_parabolics}
Suppose that $S$ is an incompressible, $\partial$--incompressible surface
which represents boundary slope $\a$.

Then there is a incompressible, $\partial$--incompressible surface $S'$
also with boundary slope $\a$ which contains no accidental parabolics.

Furthermore, $\pi_1(S')$ is a subgroup of $\pi_1(S)$.
\end{lemma}
\proof If the original surface $S$ contains no accidental parabolics,
then we are done.  Otherwise, it is shown in \cite{T} Proposition  8.11.1, that
there are two finite collections of disjoint simple closed curves $P^\pm$
so that any loop which is an accidental parabolic is freely homotopic to (a
power of) some leaf of one of these laminations.

Let $\xi$ be a leaf in one of these laminations.  We claim that the other
end of the free homotopy in $\partial M$ must be a curve which is parallel
to a (power of) $\a$.

The reason is this: The free homotopy defines a map of an annulus $f \co  A
\rightarrow M$, one of whose boundary components $f(\partial_\xi A)$ lies
on $S$ and is parallel to $\xi$ and one of whose boundary components lies
on $\partial M$.  Put $f(A)$ in general position with respect to $S$. 
Using a very small move, lift $f(\partial_\xi A)$ so that it lies slightly
away from $S$.

Since $f(A)$ and $S$ are in general position, we see $f^{-1}(S)$ as a
collection of arcs and circles in $A$.  Since $f(\partial_\xi A) \cap S$ is
empty, every arc must run from the boundary component of $A$ which is
mapped to $\partial M$ back to this boundary component.

If such an arc of intersection is essential when mapped into $S$, then it
give rise to a mapped in boundary compression and hence a compression of
$S$, since $S$ is not an annulus.  This contradicts our assumption on $S$. 
Otherwise, then the arc of intersection can be removed by homotopy.

It follows that we may suppose that there are no arcs of intersection in
$f(A) \cap S$.  However this implies that the boundary component of $f(A)$
which lies on $\partial M$ can be homotoped so as to be disjoint from
$\partial S$, proving the claim.

Consider an annulus mapped into $M$, one end mapping to $\xi$ and the other
end mapped into $\partial M$.  Put this annulus in general position with
respect to $S$; we see easily that there is an essential annulus with one
boundary component (now possibly not embedded) in $S$ and the other
boundary component a power of $\a$, so that the interior of the annulus
does not map into $S$.  Thurston's observation implies that the end which
lies on $S$ must be a power of some simple loop.  Call this simple loop
$\xi'$, say.

Let $M_{cut}$ be $M - N(S)$, where $N(S)$ is some open regular
neighbourhood of $S$.  This is an irreducible $3$--manifold with a
boundary component of positive genus, so that $M_{cut}$ is Haken.  The
observation of the previous paragraph shows that $M_{cut}$ contains an
essential annulus.

Following Jaco (see \cite{Ja}) we define a {\em Haken manifold pair}
$(X,Y)$ to be a Haken manifold $X$ together with some incompressible,
possibly disconnected $2$--manifold $ Y \subset \partial M$.  In our
setting, take $Y_1$ to be a neighbourhood of the simple closed curve $\xi'$
and $Y_2$ to be a neighbourhood of $\a$; the annulus of the above paragraph
provides an essential map of pairs $$(S^1 \times I , S^1 \times \partial I)
\longrightarrow (M_{cut} , Y_1 \cup Y_2).$$ Now Theorem VIII.13 of \cite{Ja}
implies that there is an embedding of an annulus with one end in $Y_1$ and
the other end in $Y_2$.

We now return to $M$ and use this embedded annulus to form a new embedded
surface $S_1$ by removing a neighbourhood of $A \cap S$ from $S$ and
replacing with two copies of $A$.  This surface could be disconnected, in
which case choose some component and rename this as $S_1$.  The new surface
continues to be incompressible (and therefore $\partial$--incompressible
since $S_1$ cannot be an annulus).

We may repeat this process now with the surface $S_1$.  However this
procedure must eventually terminate in a surface $S'$ which therefore
contains no accidental parabolics.  \endproof 

 We recall that
a $\pi_1$--injective surface is a {\em virtual fibre} of $M$, if there is
some finite sheeted covering of $M$, $\tilde{M}$, to which $S$ lifts and is
isotopic to the fibre of a fibration of $\tilde{M}$.  Although we do not
actually use this fact in the sequel, we note that we may deduce:
\begin{corollary}
Suppose that $S$ is a surface of maximal Euler characteristic representing 
boundary slope $\a$.  Then either $S$ is a virtual fibre of $M$, or it is 
quasifuchsian.
\end{corollary}
\proof Such a surface cannot contain an accidental parabolic, for
their removal constructs a surface of the same boundary slope of increased
 Euler characteristic.  If the surface is geometrically
finite, then as observed above, this implies that the surface is
quasifuchsian.

If the surface is geometrically infinite, this implies that it is a virtual
fibre.  This is an argument originally due to Thurston.  The ingredients
are contained in \cite{CEG} Theorem 5.2.18 and \cite{B}. \endproof

\medskip 
{\bf Remark}\qua This, together with Gabai's proof of
Property R, gives an alternative proof of a result due to Fenley \cite{F},
who proved that if $S$ is the minimal genus Seifert surface of a nonfibred
hyperbolic knot in $S^3$, then $S$ is quasifuchsian.
\medskip

 We now recall some of the results and terminology of \cite{CG}:
Suppose that we are given a curve of representations of the
fundamental group of a complete hyperbolic $3$--manifold which
contains a faithful representation---for the purposes of this paper it
always suffices to take the component containing the complete
structure.  Then we can associate to an ideal point of this curve an
action of the group on a simplicial tree.  Surfaces are constructed by
a transversality argument as subgroups of edge stabilisers.  Such an
action or surface we say is {\em associated to an ideal point}.  For
our purpose, we need only note that \cite{CG} Proposition 1.2.7 shows
given an action associated to an ideal point, a nontrivial normal
subgroup cannot fix any point of this simplicial tree.  We show:
\begin{lemma}
\label{qf_representative}
Suppose that $S$ is an incompressible, $\partial$--incompressible surface
which is associated to some ideal point and that $S$ represents boundary
slope $\a$.

Then there is a incompressible, $\partial$--incompressible quasifuchsian
surface $S'$ also with boundary slope $\a$.
\end{lemma}
\proof By hypothesis, there is a curve of representations which
yields the surface $S$ and tree $T$ via some ideal point.  We note that the
modifications in removing annuli in the proof of Lemma \ref{no_parabolics}
yield a incompressible, boundary incompressible surface $S'$ with the
property that $\pi_1(S') \leq \pi_1(S)$, so that $\pi_1(S')$ also lies in
an edge stabiliser of $T$.  We claim that this implies that the surface
group $\pi_1(S')$ cannot be geometrically infinite.

For if it were, then as above, we deduce that $S'$ is a virtual fibre of
$M$, so that there is a finite sheeted covering $M'$ of $M$ with the
property that $\pi_1(S')$ is a nontrivial normal subgroup in $\pi_1(M')$. 
However this is impossible, since restriction gives an action of
$\pi_1(M')$ on the tree $T$ and this contradicts Proposition 1.2.7 of
\cite{CG}.

We deduce that $\pi_1(S')$ is geometrically finite, it contains no
accidental parabolics, so that it is quasifuchsian, as was claimed.  \endproof

\begin{lemma}
\label{striplemma}
Suppose that $S$ is a quasi-fuchsian surface with boundary slope $\a$.

Then there is a $K$ so that if $\g$ is any simple closed curve in $\partial
M$ with $\Delta(\a, \g) > K$, we have $$\pi_1(S) \cap
\g\cdot\pi_1(S)\cdot\g^{-1} = \langle \a\rangle $$
\end{lemma}
\proof Since the image of $\pi_1(S)$ is quasifuchsian its limit set
$\Lambda(\pi_1(S))$ is a quasi-circle.  Choose basepoints and identify the
torus subgroup with $\langle \a,\b\rangle  \subset PSL(2,{\bf C})$ so that this
stabilises $\infty$ and acts on the complex plane in the upper half space
model.  The quotient $T = {\bf C}/\langle \a,\b\rangle $ is a torus and the
quasi-circle descends to a closed curve $\Lambda(\pi_1(S))/\langle \alpha\rangle  = C$
which is a homotopic to a closed geodesic on this torus.  Fix some homotopy
which moves $C$ to a Euclidean geodesic, this homotopy is covered by a
homotopy which carries $C$ to a straight line $s_C$.  Since the homotopy is
the image of a compact set, this shows that the set $\Lambda(\pi_1(S))$
lies within some fixed distance $\ell$ of $s_C$ in the Euclidean metric on
${\bf C}$.

Choose $K$ so large that translation by $\b^K$ moves the line $s_C$ a
distance of $1000\ell$ away from itself.
If we choose an element $\g$ as in the statement of Theorem \ref{main},
(ie, $\g = \a^r\beta^n$, where $|n| > K$), it follows that the limit set
of the subgroup $\pi_1(S)$ and the limit set of the subgroup
$\g\cdot\pi_1(S)\cdot\g^{-1}$ contain only the point at infinity in common. 
It follows that the only possibility for elements in $\pi_1(S) \cap
\g\cdot\pi_1(S)\cdot\g^{-1}$ lie in the parabolic subgroup fixing $\infty$
since such elements must stabilise both limit sets; thus the intersection
contains $\langle \a\rangle $, but can be no larger.  \endproof 

 Notice
that the proof actually constructs a straight strip in ${\bf C}$ which
contains the limit set of $\pi_1(S)$.  We also observe that the $K$ of the
lemma need not yet be the $K$ of Theorem \ref{main}, we may need to enlarge
it further.

To be specific, let us henceforth suppose that the quasifuchsian surface
$S$ contains two boundary components, the case of a single boundary
component having been dealt with in \cite{CL} and the case of more boundary
components being entirely analogous.

\subsection{The complex}

We now begin our construction of a certain complex ${\cal C}$.

The surface $S$ is finite area and quasifuchsian.  We recall that
$\Hull(\Lambda(\pi_1(S))$ is the intersection of all the hyperbolic
halfspaces which contain $\Lambda(\pi_1(S)$; this is a convex $\pi_1(S)$
invariant set.  In the degenerate case that the surface is totally 
geodesic, we adopt the convention that $\Hull(\Lambda(\pi_1(S))$ is
defined to be small $\epsilon$-neighbourhood of the intersection
of the hyperbolic halfspaces containing $\Lambda(\pi_1(S))$.

The fact that $S$ is quasifuchsian means that
$\Hull(\Lambda(\pi_1(S))/\pi_1(S)$ is a finite volume hyperbolic manifold
with cusps homeomorphic to $S \times I$.  We may homotope $S$ so that $\tilde{S}$ 
lies inside its hull, for example by homotoping it just inside one of the
pleated surface boundaries. By further homotopy, we arrange that  $S$ meets
the boundary of a cusp in simple closed curves. 

From this it follows that we can bound the ``thickness" of the hull, that
is to say, there is some constant, $c_1$, so that every point on one of the
surfaces of the convex hull of $\tilde{S}$ is a distance at most $c_1$ from
the other surface in the convex hull and in particular, every point in the
hull is within distance $c_1$ of some point of $\tilde{S}$.  We may as well
assume that $c_1$ is fairly large, at least $10$ say.

Fix a horoball neighbourhood in $M$ of the cusp, $N'$, so small that the
smallest distance between two preimages of $N'$ in the universal covering
is very large compared to $1000c_1$.  We can assume that $N'$ was chosen
small enough so that only the thin part of $S$ coming from the cusps enters
$N$.  Now chose a horoball $N$ so far inside $N'$ that the hyperbolic
distance between $\partial N'$ and $\partial N$ is greater than $1000c_1$.

This now guarantees that in the universal covering, every preimage of $N$
either is centred at some limit point of $\tilde{S}$ or is distance greater
than $1000c_1$ from $\tilde{S}$.  Using our thickness estimate, we see that
every preimage of $N$ is either centred at a limit point of $\tilde{S}$ or
it is at a distance much greater than $750c_1$ from the convex hull of
$\tilde{S}$.

Define $M^-$ to be the complete hyperbolic manifold $M$ with the interior
of $N$ excised.  This is a compact manifold with a single torus boundary
component.  We have arranged the surface $S$ so that the surface $S^-$ has
two boundary components on this torus.  We will base all fundamental groups
at some point $p$ on one of these boundary components.

The universal covering of $M^-$ embeds into the universal covering of $M$
in the obvious way.  We use the upper halfspace model and let $\tilde{N}$
be the horoball preimage of $N$ centred at infinity.  Fix some preimage of
$S$, $\tilde{S}_1$ which passes through infinity.  As in the proof of Lemma
\ref{striplemma}, we may find a pair of vertical planes which define a
three dimensional strip which meets ${\bf C}$ in a strip containing the
limit set of $\tilde{S}_1$.

In the upper halfspace model all horoballs not centred at $\infty$ have
some uniformly bounded Euclidean size, so by moving these planes apart if
necessary, we may suppose that the three dimensional strip contains the
intersection of $\tilde{S_1}$ with the horosphere corresponding to
$\tilde{N}$ as well as all the translates of $\tilde{N}$ which meet
$\tilde{S}_1$.  Denote the three dimensional strip this produces by
$\Sigma_1$.  The intersection of $\Sigma_1$ with the $\tilde{N}$ horosphere
will be denoted $\sigma_1$.

Notice that the set $\Sigma_1$ is bounded by two totally geodesic planes,
so that it is hyperbolically convex, moreover it contains the limit set of
$\tilde{S}_1$ by construction, so that it contains the convex hull of this
limit set.  The set $\tilde{S}_1$ meets the horosphere $\partial\tilde{N}$
in a line which covers one of the boundary components of $S^-$, let us
denote this component by $\partial_1S^-$.  This situation is depicted in
Figure 1.
\begin{figure}[ht!]
\cl{\relabelbox\small
\epsfysize=200pt
\epsfbox{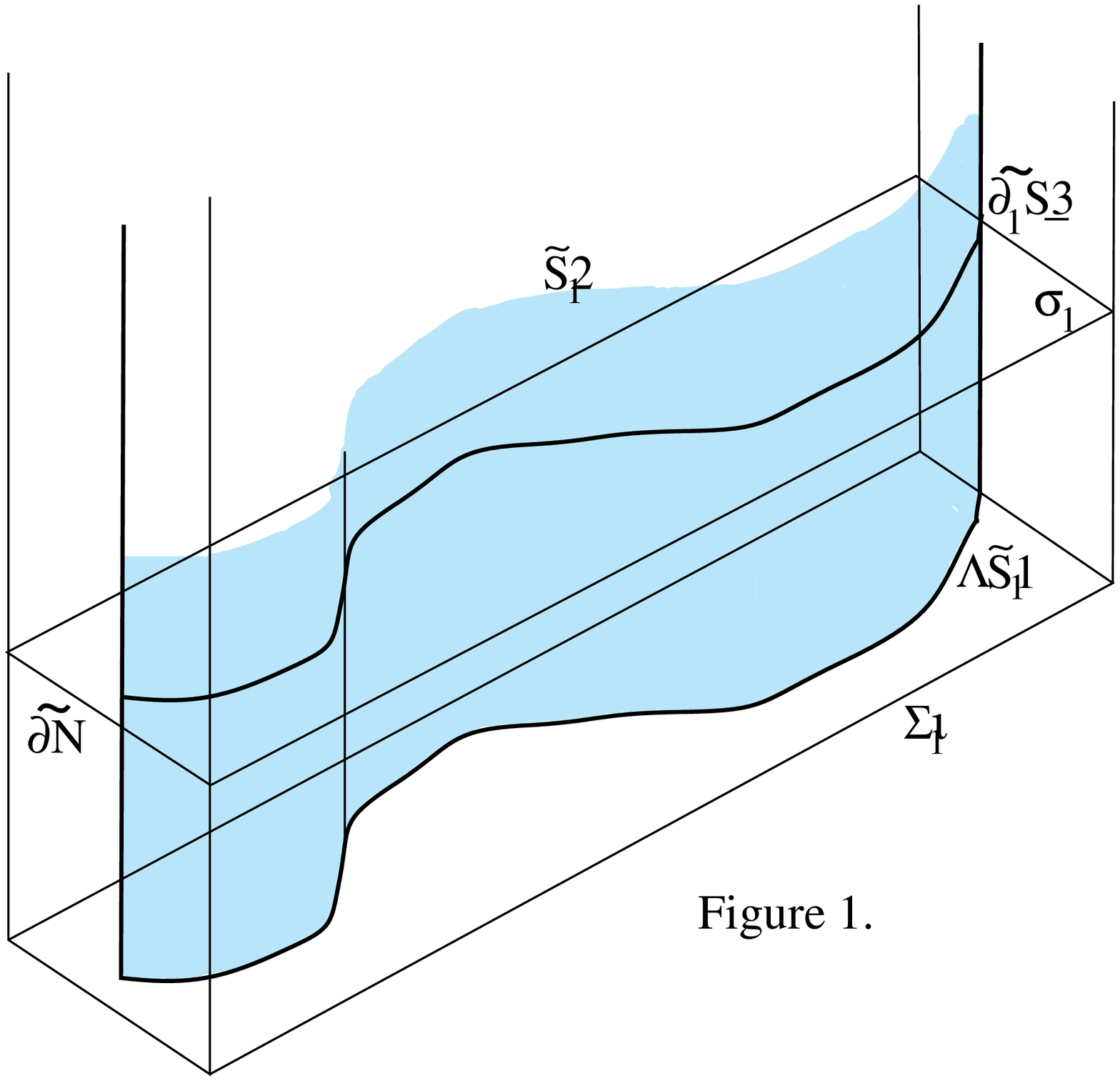}
\relabela <-4pt,0pt> {S1}{$\Lambda\widetilde S_1$}
\relabela <0pt,3pt> {S2}{$\widetilde S_1$}
\relabel {Si}{$\Sigma_1$}
\relabel {S3}{$\widetilde{\partial_1 S_-}$}
\relabela <0pt,-3pt> {N}{$\widetilde{\partial N}$}
\relabel {s}{$\sigma_1$}
\endrelabelbox}
\nocolon\caption{}
\end{figure}

Now we fix a $\pi_1(M)$--translate of $\tilde{S}$ which passes through
infinity and meets $\partial\tilde{N}$ in a line which covers the other
boundary component $\partial_2S^-$.  Denote this translate by
$\tilde{S}_2$.  We then perform the construction of the above paragraph
with $\tilde{S}_2$ to form the sets $\Sigma_2$ and $\sigma_2$.  Note that
$\sigma_1$ and $\sigma_2$ are parallel strips, so by applying some
preliminary covering translation if necessary, we may assume that
$\tilde{S}_2$ was chosen so that the distance between $\sigma_1$ and
$\sigma_2$ is very large in both the hyperbolic metric and in the Euclidean
metric of the horosphere. (See Figure 2A.)

The proof of \ref{striplemma} together with the fact that the translation
coming from $\alpha$ stabilises the strips implies that we may choose a $K$
so that if $\g$ is in the stabiliser of $\infty$ and $\Delta(\a,\g) > K$,
then the distance between $\sigma_1 \cup \sigma_2$ and $\g(\sigma_1 \cup
\sigma_2)$ is very much larger than the distance between $\sigma_1$ and
$\sigma_2$ in both the hyperbolic metric and in the Euclidean metric of the
horosphere $\partial\tilde{N}$.  (In particular, this distance is very
large.)  Since this distance is very large, we may as well assume at the
same time that if $\Delta(\a,\g) > K$, then the length of $\g$ in the
Euclidean metric on $\partial\tilde{N}$ is much larger than $2\pi$.

Suppose that we wish to do Dehn filling corresponding to some curve $\g
\subset \partial M^-$ which satisfies the requirement that $\Delta(\a,\g) >
K$, with the choice of $K$ in the above paragraph.

Consider the subset of hyperbolic $3$--space defined as follows.  Begin
with the two convex hulls coming from $\tilde{S}_1$ and $\tilde{S}_2$. 
Adjoin all the horoballs which cover $N$ and which are incident on
$\tilde{S}_1$ or $\tilde{S}_2$.  Note that this configuration is invariant
under the element $\a$.  Now add in all $\g$ translates in the
configuration so that it becomes invariant under the subgroup $\langle \a,\g\rangle $. 
Denote this set by ${\cal C}_0$.

Suppose we have inductively constructed the complex ${\cal C}_n$; this
consists of certain translates of the convex hulls of $\tilde{S}_1$ and
$\tilde{S}_2$ and certain translates of horoballs which cover $N$. 
Inductively, we assume that the entire complex is invariant under the group
$\langle \a,\g\rangle $.

Some of these horoballs of ${\cal C}_n$ are centred at points which meet
translates of the limit sets of both $\tilde{S}_1$ and $\tilde{S}_2$; that
is to say, the surface $S$ has two ends (corresponding to the fact that
we're assuming that $S$ has two boundary components) and both of these ends
appear in this type of horoball.  (By way of example, $\infty$ is the only
such horoball in ${\cal C}_0$.)  Inductively we may assume that for such
horoballs, the collection of surfaces which is incident to the horoball is
invariant under the orbit of the conjugate of the subgroup $\langle \a,\g\rangle $ which
stabilises the horoball.

However some of the horoballs of ${\cal C}_n$ only meet translates of the
preimage of one end of $S$.  It is along these horoballs that we enlarge
the complex: To be specific, we will denote such a surface by
$\tilde{S}(j)$, and the horoball by $\tilde{N}(j)$, assuming it is centred
at the complex number $j \in S^2_\infty$.

There are now two cases.  (In general, the number of cases is the number of
components of $\partial S$.)

The first case is that $\tilde{S}(j)$ meets $\tilde{N}(j)$ in some closed
curve which covers $\partial_1 S^-$.  In this case we may actually choose
an element of $\pi_1(M)$ which maps $\infty$ to $j$ while mapping
$\tilde{S}_1$ to $\tilde{S}(j)$.  The ambiguity in such an element is easily
seen to come from premultiplication by any element of $\langle \a\rangle $.

The second case is that $\tilde{S}(j)$ meets $\tilde{N}(j)$ in some closed
curve which covers $\partial_2 S^-$.  In this case we may choose an element
of $\pi_1(M)$ which maps $\infty$ to $j$ while mapping $\tilde{S}_2$ to
$\tilde{S}(j)$.  Again, the only ambiguity in such an element comes from
premultiplication by an element of $\langle \a\rangle $.

In either case, denote some choice of such an element by $\mu_j$.  We form
${\cal C}_{n+1}$ by adding to ${\cal C}_n$ all the complexes $\mu_j\cdot
{\cal C}_0$ as we run over all the relevant horoballs.  Since ${\cal C}_0$
is $\langle \a,\g\rangle $--invariant , this is independent of the $\mu$--choices.  One
sees easily from the construction that ${\cal C}_{n+1}$ satisfies the
inductive hypothesis.

Define ${\cal C}$ to be the union of all the ${\cal C}_n$'s.  A schematic
for ${\cal C}$ is shown in Figure 2A.

To expedite our analysis of ${\cal C}$, we define a graph $\cal G$ by
taking one vertex for each preimage of $N$ lying in ${\cal C}$ and one
vertex for each translate of $\tilde{S}$ in ${\cal C}$.  Edges are defined
by the obvious incidence relations.  By construction, $\cal G$ has two
types of vertex and two types of edge.  (In general, there are two types of
vertex and $|\partial S|$ types of edge.)

\begin{figure}[ht!]
\cl{\relabelbox\small
\epsfxsize=.47\hsize
\epsfbox{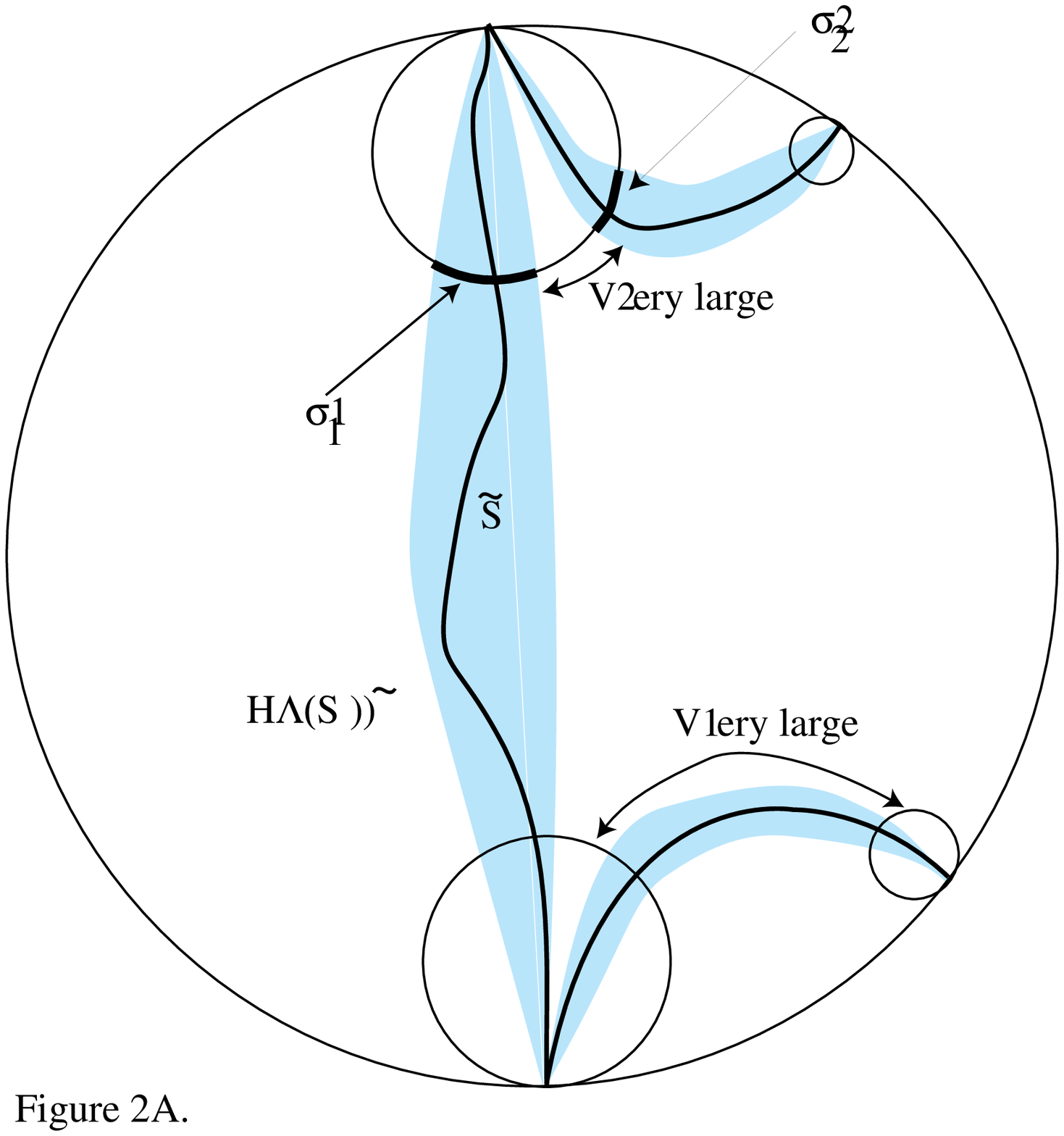}
\relabel {V1}{Very large}
\relabela <0pt,-3pt> {V2}{Very large}
\relabela <0pt,-1pt> {s1}{$\sigma_1$}
\relabela <-2pt,0pt> {s2}{$\sigma_2$}
\relabela <-15pt,0pt> {H}{Hull$(\Lambda(\widetilde S))$}
\relabel {S}{$\widetilde S$}
\endrelabelbox
\relabelbox\small
\epsfxsize=.47\hsize
\epsfbox{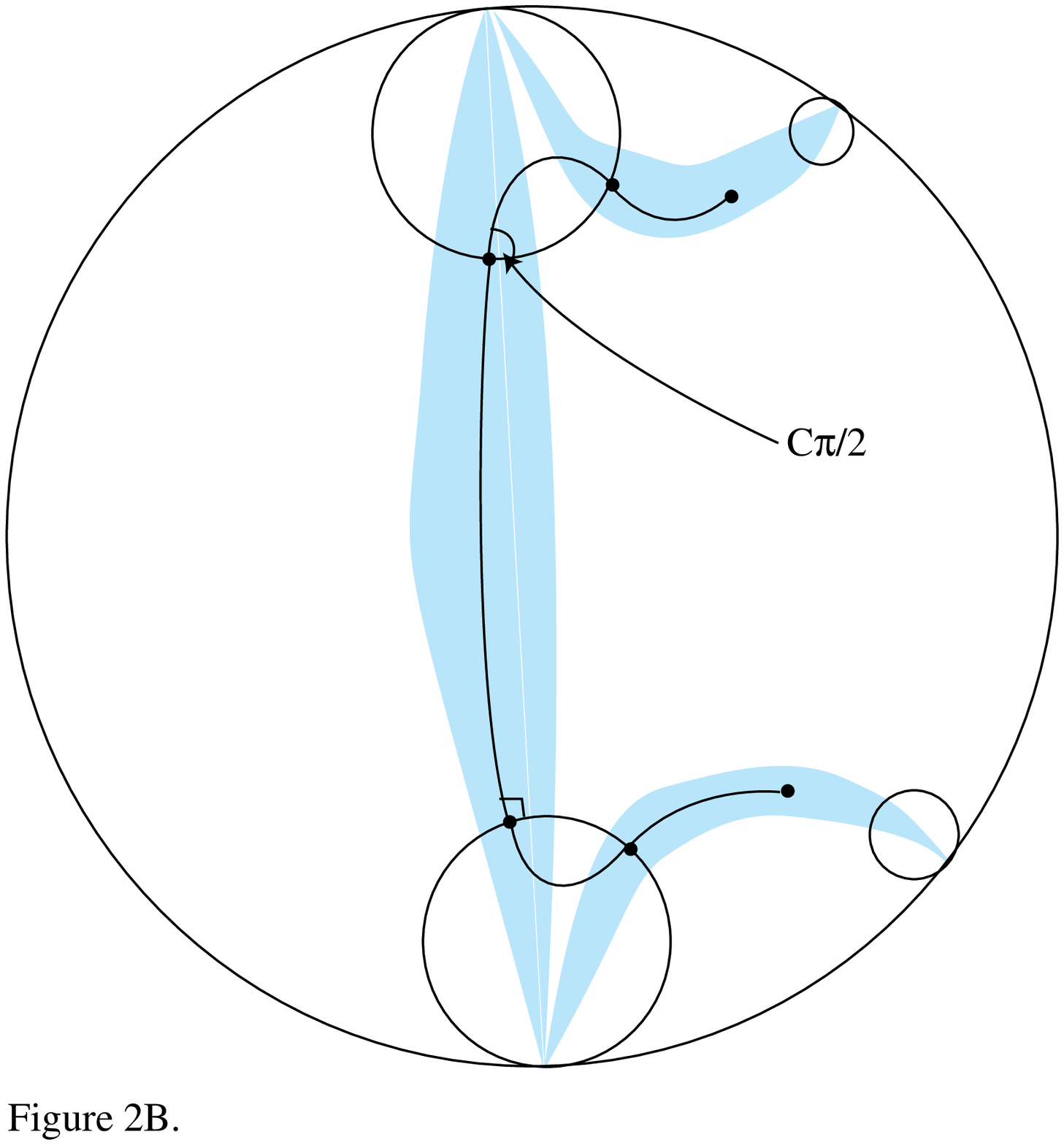}
\relabela <-15pt,-6pt> {C}{Close to $\pi/2$}
\endrelabelbox}
\nocolon\caption{A \hglue 2in Figure 2B}
\end{figure}

\begin{lemma}
\label{graph_is_a_tree}
The graph $\cal G$ is a tree.
\end{lemma}
\proof Suppose not and that there is a path in $\cal G$ which runs
between vertices without backtracking.  This corresponds to some path in
${\cal C}$ which runs from one preimage of the surface $S$ to another.  We
will replace this path by a piecewise geodesic path which will in fact be a
long quasi-geodesic.  From this it will follow that the path runs between
different lifts of $S$, so that it is not a loop in $\cal G$.

Here is the construction: Since it does not involve backtracking, such a
path consists of pieces of two types.  There are segments which run from
horoball to horoball inside a convex hull of some lift of $S$, and segments
which run inside the horoball from one lift of $S$ to another.

If we see a segment of the first type, we replace it by the common geodesic
perpendicular to the pair of horoballs.  Notice that this geodesic runs
through the centres of both horoballs, which lie in the limit set of this
preimage of $S$, so that the geodesic continues to lie inside the convex
hull.  By construction, such a geodesic makes angle $\pi/2$ with the
horosphere that it meets.

If we see a segment of the path lying inside the horoball connecting
different lifts of $S$, we replace it by the geodesic in the horoball
between endpoints of the relevant common perpendiculars constructed in the
previous paragraph.  Note that the distance between such points is
enormous, so that this geodesic makes an angle with the horosphere which is
very close to $\pi/2$.  (See Figure 2B.)

Our  choices ensure that all the geodesics we construct this way are
all very long, moreover, the angle a horoball geodesic piece makes with a
convex hull geodesic piece is very close to $\pi$.  This makes the path
quasi-geodesic, hence it has distinct endpoints and does not correspond
to a closed path in the graph $\cal G$ as was required.  \endproof

{\bf Remark}\qua As a first application of this lemma, we
note that it implies that if $j$ is the centre of some horoball in ${\cal
C}$, then there are no ``unexpected" translates of $\tilde{S}$ in ${\cal
C}$ incident on $j$.  That is to say, in the above notation, if
$\mu_j(\infty) = j$, then the first appearance of $j$ in our construction
places into ${\cal C}$ the $\mu_j$ image of the orbit $\langle \a,\g\rangle (\tilde{S}_1
\cup \tilde{S}_2)$ at $j$.  These are the only preimages of $S$ in ${\cal
C}$ incident at $j$; the reason being that any preimage added at a later
stage in the construction would give rise to a loop in the graph $\cal G$.

\begin{corollary}
\label{stabcorollary}
Suppose that $g \in \pi_1(M)$ carries some preimage of $S$ lying in ${\cal C}$
to another preimage of $S$ which lies in ${\cal C}$.

Then $g$ stabilises ${\cal C}$.
\end{corollary}
\proof We begin by noting that ${\cal C}$ has the following property:
Suppose that $j$ is the centre of some horoball in ${\cal C}$,
$\tilde{N}(j)$ say.  This means that there is some translate of
$\tilde{S}$, $\tilde{S}(j)$ say, in the complex ${\cal C}$ whose limit set
contains $j$.  We claim that knowing any translate determines the entire
complex ${\cal C}$.

The reason for the claim is this: The translates of $\tilde{S}$ which lie
in ${\cal C}$ and have limit point at $j$ are constructed as the images of
the $\langle \a,\g\rangle $ orbit of $\tilde{S}_1 \cup \tilde{S}_2$ under some element
$\mu_j$, which in particular satisfies $\mu_j(\infty) = j$.  By the remark
following \ref{graph_is_a_tree}, this is exactly the collection of
translates in ${\cal C}$ at $j$.

Since $\tilde{S}(j)$ lies in ${\cal C}$, it follows that we can find an
element $\xi \in \langle \a,\g\rangle $ so that $\tilde{S}(j) = \mu_j\xi\tilde{S}_*$,
where $\tilde{S}_*$ denotes one of the reference surfaces, ie, either
$\tilde{S}_1$ or $\tilde{S}_2$.  (Of course, which one is determined by
$\tilde{S}(j) \cap \partial \tilde{N}(j)$ and which component of $\partial
S^-$ that this covers.)

Suppose that $g_j$ is any element of $\pi_1(M)$ which throws $\infty$ to
$j$ and the relevant reference surface to $\tilde{S}(j)$.  The element
$g_j^{-1}\mu_j\xi$ stabilises $\infty$ and the reference surface, so it is
a power of $\a$, that is to say $g_j = \mu_j\xi\a^t$ for some integer $t$. 
We deduce that we can reconstruct the translates of $\tilde{S}$ in ${\cal
C}$ at $j$ as the $g_j$ image of the orbit $\langle \a,\g\rangle \cdot(\tilde{S}_1 \cup
\tilde{S}_2)$.

It follows then that we can unambiguously reconstruct the complex ${\cal
C}$ outwards using the new set of horoballs that this created, since each
such horoball meets a translate of $\tilde{S}$ and we may apply the same
argument.  This proves our claim.

The corollary now follows: If $g$ is some element of $\pi_1(M)$ throwing
$\tilde{S}(i)$ to $\tilde{S}(j)$ and $i$ to $j$, then in the notation
established in this proof, we may choose $g_j$ to construct the surfaces in
${\cal C}$ incident to $j$ and $g_j\circ g$ to construct all the surfaces
in ${\cal C}$ incident to $i$.  From this it follows that $g$ maps the
surfaces in ${\cal C}$ incident to $i$ to the surfaces in ${\cal C}$
incident to $j$.  The result follows by building ${\cal C}$ outwards as
described above.  \endproof

 We will need to know that $\Hull({\cal C})$ has the following property. 
\begin{theorem}
\label{Nproperty}
The set $\Hull({\cal C})$ has the property that any $\pi_1(M)$ translate of
$\tilde{N}$ is either entirely contained inside $\Hull({\cal C})$ or is
disjoint from it.
\end{theorem}
The analysis which achieves this is technical and will be deferred to 
section \ref{technicalsection}

\subsection{The subgroup $\Gamma$}
Define $\Gamma$ to be the stabiliser in $\pi_1(M)$ of the complex ${\cal
C}$.  We will show that $\Gamma$ is isomorphic to the fundamental group of
a certain $3$--manifold.  There is a minor difference between the cases
that $S$ is the boundary of a twisted $I$-bundle or it is  not and so we shall assume
henceforth: $S$ is not the boundary of a twisted $I$--bundle neighbourhood
of a nonorientable surface embedded in $M$.

The minor changes which need to be made in this situation are explained in
section \ref{twistedIbundles}.  In fact the only case which really needs to be
singled out is the case when $S$ is the boundary of a twisted $I$--bundle
neighbourhood of a nonorientable surface {\em with one boundary component}.

Define this $3$--manifold as follows: Take $T^2 \times [0, \infty)$ and
glue on a manifold homeomorphic to $S^- \times I$ by identifying $\partial
S^- \times I$ with two annuli parallel to $\a$ contained in $T^2 \times 0$. 
Denote this manifold by $A$.  Of course, a manifold homeomorphic to $A$ is
already embedded inside $M^-$, however the group $\Gamma$ corresponds to an
immersion of $A$ which intuitively comes by attaching one boundary
component of $\partial S^-$ to $\partial M^-$, then spinning the other
boundary component around a very long annulus parallel to $\partial M^-$
before attaching to $\partial M^-$.

The manifold $A$ is easily seen to be irreducible and its boundary contains
a component of positive genus, so that $A$ is Haken.  We will show that
$\Gamma \cong \pi_1(A)$.  To this end we need:
\begin{lemma}
The complex ${\cal C}$ is simply connected.
\end{lemma}
\proof We note that ${\cal C}$ is constructed as horoballs glued to
copies of $\widetilde{S^-} \times I$ along thin neighbourhoods of copies of
the real line which cover boundary components of $S^-$.  The result then
follows from the Seifert--Van Kampen theorem, together with Lemma
\ref{graph_is_a_tree}.  \endproof

 We may now prove:
\begin{theorem}

 The complex ${\cal C}/\Gamma$ is homeomorphic to $A$, in particular,
 $\Gamma \cong \pi_1(A)$.

\end{theorem}
\proof The group $\Gamma$ acts on ${\cal C}$ and Corollary
\ref{stabcorollary} shows that
\begin{enumerate}
\item $\Gamma$ acts transitively on the copies of $\tilde{S}$ lying in
${\cal C}$.
\item $stab(\tilde{S})$ is a subgroup of $\Gamma$.
\end{enumerate}
Taking together $1.$ and $2.$ we see that $\Gamma$ acts transitively on the
horoballs of ${\cal C}$.  Moreover, since we assumed that $S$ was not the
boundary of an $I$--bundle, it follows from results in Chapter $10$ of
\cite{He}, that $stab(\tilde{S}) \cong \pi_1(S)$.  Given this we see that
${\cal C}/\Gamma$ is homeomorphic to $A$.  Since the complex ${\cal C}$ is
simply connected, the result follows.  \endproof 

{\bf Remark}\qua The proof shows that this isomorphism identifies the
cusp subgroup of $\Gamma$ with the torus subgroup of $\partial A$ in
the obvious way.

\subsection{Constructing the surface group}
We begin with a simple lemma (see also \cite{CL}, Proposition 2.1):
\begin{lemma}
\label{containssurfacegroup}
Suppose that $F$ is an orientable surface with $k > 1$ boundary components. 
Fix some simple closed curve $C$ in the boundary of a solid torus $T$ which
does not bound a disc in $T$.  Form a $3$--manifold $P$ by identifying all
the annuli $\partial_i F \times I \subset F \times I$ with disjoint annuli
which are all parallel to neighbourhoods of $C$ in $\partial T$.

Then $\pi_1(P)$ contains the fundamental group of a closed orientable
surface.
\end{lemma}
\proof By passing to a covering if necessary, it is easy to see that
we may assume that $C$ meets the disc in the solid torus exactly once
transversally.

We form a covering space of $P$ in the following way: Take $k$ copies of
the solid torus and $k$ copies of $F \times I$.  Any way that we glue these
objects together subject to the obvious restrictions coming from the way $F
\times I$ is glued on to $T$ will be a covering space of $P$.

Fix a copy $(F \times I)_1$ and glue up in any way (subject to being
compatible with being a covering) so that each of the boundary components
of $\partial F$ appears on a different torus.  Now glue on a second copy
$(F \times I)_2$ similarly.  We claim that this $3$--manifold already
contains a closed surface group, which will imply the result, since the
remaining part of the covering is constructed by forming HNN construction
along curves which are nontrivial in both groups.

We see the claim by noting that all we have done abstractly is take two
copies of $F \times I $ and identified some component of $\partial F \times
I$ of one copy with some component of $\partial F \times I$ in the other. 
This manifold has incompressible (though possibly nonorientable) boundary.
In particular, it contains a closed surface group, proving the
result.  \endproof

 Consider the manifold $X = \Hull({\cal
C})/\Gamma$.  Of course, $\Hull({\cal C})$ is convex and therefore
contractible, so $X$ is a $K(\Gamma,1)$, moreover the isomorphism of
$\pi_1(A)$ with $\Gamma$ is the obvious isomorphism between the boundary
torus of $A$ and the cusp group $\langle \a,\g\rangle $ of $\pi_1(X)$.  For future reference we
note that $X$ is a hyperbolic manifold with convex boundary.

The restriction of the covering map ${\bf H}^3/\Gamma \rightarrow M$ 
 is a map $f \co  X \rightarrow M$ which is a local isometry.

By Theorem \ref{Nproperty}, every translate of $\tilde{N}$ is a horoball
in ${\bf H}^3$ which is either contained in $\Hull({\cal C})$ or is
disjoint from it.  It follows that we may excise $\tilde{N}/\langle \a,\g\rangle $ from
$X$ and $N$ from $M$ and define a new map $$f^- \co  X^- \longrightarrow M^-$$
which continues to be a local isometry.  This map induces a covering of
boundary tori of degree given by the index of $\langle \a, \g \rangle  $ in
$\pi_1(\partial M^-)$.  Since the curve $\g$ is mapped 1--1 by this
covering, we may extend to a map of the surgered manifolds $$f_\g \co  X(\g)
\longrightarrow M(\g)$$ 

Recall that by construction $\gamma$ has length more than $2\pi$ on
$\tilde{N}$, so that using the $2\pi$--theorem \cite{BH}, we may put
negatively curved metrics on both spaces and arrange that the map $f_\g$
continues to be a local isometry.  We now may prove:
\begin{theorem}
\label{injection}
The map $f_{\g*} \co  \pi_1(X(\g)) \longrightarrow \pi_1(M(\g))$ is injective.
\end{theorem}
\proof We begin by observing that although the manifold $X^-$ is not
convex, the manifold $X(\g)$ with the extended metric is convex, since
$\partial X(\g)$ is convex.

Suppose then that the theorem were false and that we could find some
element in the kernel of $f_{\g*}$.  Since $X(\g)$ is convex, any such
element is freely homotopic to a geodesic $\xi$ in $X(\g)$.

The loop $f_\g(\xi)$ is a geodesic in $M(\g)$, since $f_\g$ is a local
isometry.  This is a contradiction; if the loop $f_\g(\xi)$ were
nullhomotopic, it would lift to be a nullhomotopic geodesic in
$\widetilde{M(\g)}$, a negatively curved complete simply connected manifold
and this contradicts the theorem of Hadamard-Cartan.  (\cite{dC} Theorem
3.1).  \endproof

 Theorem \ref{main} now follows:
\begin{theorem}
The manifold $M(\g)$ contains the fundamental group of a closed orientable
surface
\end{theorem}
\proof The manifold $\pi_1(A(\g)) \cong \pi_1(X(\g))$ contains the
fundamental group of a closed surface by Lemma \ref{containssurfacegroup}
and this group injects into $\pi_1(M(\g))$ by Theorem \ref{injection}. 
\endproof

{\bf Remarks}\qua The reason for using $\Hull({\cal C})$ rather than, for
example $\Hull(\Lambda(\Gamma))$ is that there seems to be no a
fortiori control over the size of the horoball which embeds into this
latter set.  Such control is needed in order to apply the $2\pi$ theorem.

It also seems worth clarifying why we do not work with the map 
${\bf H}^3/\Gamma \rightarrow M$ directly. The reason is that in this 
setting, to define a map into $M^-$  it is necessary to excise the full
 preimage in ${\bf H}^3$ of the horoball $\tilde{N}$. The group $\Gamma$ stabilises
no parabolic fixed points of $\pi_1(M)$ which lie outside $\Lambda(\Gamma)$,
so there is no way to extend the surgery over those horospheres---and the 
resulting manifold cannot be made to have convex boundary.

\subsection{The case that $S$ is the boundary of a twisted $I$--bundle}
\label{twistedIbundles}
A variation on the construction above is necessary in the case that that
$S$ is the twisted $I$--bundle neighbourhood of a nonorientable surface
since the stabiliser of $\tilde{S}$ is a ${\bf Z}_2$ extension of the
surface group $\pi_1(S)$. It follows that the manifold $A$ as described
above is not the correct model for the fundamental group.

The above proof is easily modified if the surface core of the $I$-bundle has
at least two boundary components.  One performs an analogous construction
by gluing a twisted $I$-bundle neighbourhood onto a solid torus.  With
this proviso, we continue to have ${\cal C}/\Gamma \cong A$ and the
manifold $A$ continues to contain a surface group.

However, if the nonorientable surface core of the $I$--bundle has only a
single boundary component, then $X$ is double covered by a handlebody.  We
now sketch how to modify the above proof to deal with this case, and in
fact the more general case that $S$ is the boundary of the twisted
$I$--bundle over a nonorientable surface $F$ and that $F$ has $k \geq 1$
boundary components.

The construction of ${\cal C}$ is a mild variation on the construction
described above.  Take $k$ ``white" copies of the universal covering of
$F$, which pass through $\infty$ and meet the boundary of a preimage
horoball $\tilde{N}$ (chosen to be small as above) in lines which cover the
$k$ boundary components of $F^-$.  As usual arrange that these $k$ copies
are very far apart.  Similarly, take $k$ ``black" copies and assume that
the black and white complexes are very distant from each other.

Now choose a $K$ so that if $\Delta(\a,\g) > K$, then $\g $ maps this
complex of $2k$ surfaces a very long distance away from itself.

The rest of ${\cal C}$ is constructed much as before save only that we do
not allow the identification of white and black surfaces.  The group
$\Gamma$ is defined to be the stabiliser of the {\em coloured} complex. 
The proof that the associated graph is a tree is identical.  However the
complex ${\cal C}/\Gamma$ is a little different and for this we need to
modify $A$.  The reason for this difference is that the stabiliser of the
vertex corresponding to $\tilde{F}$ is no longer $\pi_1(S)$, but the
two-fold extension $\pi_1(F)$.

This having been noted, we define the complex $A$ to be $T^2 \times I$ with
{\em two} copies of the twisted $I$--bundle over $F$ attached in the
obvious way.  That this complex contains a surface group now can be proved
as in \ref{containssurfacegroup}.  The rest of the proof is identical.

\medskip 
{\bf Remark}\qua Though this seems artificial, we note that
in the case of a surface with a single boundary component (eg a Seifert
surface, see \cite{CL}) some sort of operation of this sort is necessary;
one cannot tube a surface with a single boundary component to itself.

\subsection{The proof of Theorem \ref{Nproperty}}
\label{technicalsection}
This section is devoted to the proof of the technical result \ref{Nproperty},
which for the reader's convenience we restate here:
\begin{theorem}
The set $\Hull({\cal C})$ has the property that any $\pi_1(M)$ translate of
$\tilde{N}$ is either entirely contained inside $\Hull({\cal C})$ or is
disjoint from it.
\end{theorem}
We recall that the {\em thin triangles constant}, $\Delta$, for hyperbolic
space, is a constant with the property that every point on one side of a
geodesic triangle is within a distance $\Delta$ of some point on the union
of the other two sides.  In fact $\Delta < 3.$

\begin{lemma}
\label{firstestimate}
There is a constant $K_1 =1+2\Delta\leq 10$ with the following property. 
Suppose that $a$ and $b$ are two points in ${\cal C}$.  Then every point on
the geodesic segment connecting $a$ to $b$ is within a distance $K_1$ of
${\cal C}$.
\end{lemma}
\proof Throughout this proof we will refer to a translate of the
chosen horoball $\tilde{N}$ simply as a {\em horoball,} and to a translate
of ${\cal C}\cap \Hull(\tilde{S})$ as a {\em thick surface.}

In \ref{graph_is_a_tree} we constructed a piecewise geodesic, $\b$ in
${\cal C}$ with endpoints $a$ and $b.$ The geodesic segments in this path
are of two types.  The first type is a geodesic segment in a thick surface. 
The second type is a geodesic segment in a horoball.  The endpoints of
these segments, other than $a$ and $b,$ are on the intersection of some
thick surface with some horoball.  A segment contained in a thick surface
is orthogonal to the boundary of any horoballs at its endpoints.  The same
is not generally true for segments in horoballs.

In \ref{graph_is_a_tree} every segment in a horoball was long (length at
least $1000$) and therefore almost orthogonal to the boundary of the
horoball at its endpoints.  In the present situation, this fails precisely
when $a$ or $b$ is in a horoball.  Thus $\b$ is a piecewise geodesic, and
all the segments, except possibly the first and last, have length at least
$1000.$ The angle between two adjacent segments is very close to $\pi$
 except possibly for the first and last angles.

Using the above properties of $\b$ it is well known that the geodesic with
endpoints $a,b$ lies within $2\Delta+1$ of $\d$.  We sketch this: Let
$[a,c]$ be the first geodesic segment and $[e,b]$ the last segment of $\b$. 
The subpath $\b_-$ of $\b$ excluding the first and last segment consists of
segments all of length at least $1000$ and with all angles between adjacent
segments within the range $\pi \pm 0.001.$ Let $[c,e]$ be the geodesic with
endpoints $c$ and $e$.  Then every point of $\b_-$ is within $1$ of $[c,e].$
Now consider the geodesic path $\delta=[a,c].[c,e].[e,b].$ Let $[a,b]$ be
the geodesic with the same endpoints as $\delta.$ Then $[a,b].\delta$ is a
geodesic quadrilateral.

By dividing this into two geodesic triangles, we see that every point on
$[a,b]$ is within $2\Delta$ of $\delta.$ Now every point of $\delta$ is
within $1$ of a point of $\beta.$ Hence every point of $[a,b]$ is within
$1+2\Delta$ of $\beta$ which is contained in $\cal C$ \endproof

\begin{lemma}
\label{secondestimate}
There is a constant $K_2=K_1+2\Delta\leq 20$ such that if $T$ is a geodesic
tetrahedron with all four vertices in ${\cal C}$ then every point of $T$ is
within a distance $K_2$ of ${\cal C}$.  \end{lemma} 

\proof We prove below that every point in the boundary of $T$ is
within a distance $K_1+\Delta$ of some point in ${\cal C}.$ Given a
point $p$ in the interior of $T,$ there is a geodesic triangle with
edges contained in the boundary of $T$ and which contains $p.$ Hence
$p$ is within a distance of $\Delta$ of some point in the boundary of
$T.$ Therefore within $K_1+2\Delta$ of some point of ${\cal C}.$

Let $F$ be a face of $T.$ Thus $F$ is a geodesic hyperbolic triangle with
vertices $a,b,c$ contained in ${\cal C}.$ Suppose that $z\ne c$ is a point
in $F.$ There is a geodesic segment, $[c,y]$ in $F$ containing $z$ and with
one endpoint $c$ and the other endpoint, $y,$ on $[a,b].$ Let $p$ be a
point of ${\cal C}$ closest to $y.$ By Lemma \ref{firstestimate},
$d(p,y)<K_1.$ Moreover, $z$ is within $\Delta$ of $[y,p]\cup [p,c].$ The
points $p,c$ are both in ${\cal C}$.  With another application of Lemma
\ref{firstestimate}, we see that every point on $[p,c]$ is within a
distance $K_1$ of ${\cal C}.$ Since $p$ is in ${\cal C}$ every point on
$[p,y]$ is within $d(p,y)<K_1$ of ${\cal C}.$ Therefore $z$ is within a
distance $\Delta$ of a point within a distance $K_1$ of ${\cal C}.$ Hence
$d(z,{\cal C})<K_1+\Delta.$ \endproof

\begin{lemma}
\label{hullobservation}
Suppose that $C_n$ is an increasing sequence of compact subsets of ${\bf
H}^3.$ The convex hull, $H,$ of $\cup C_n$ is the closure of $\cup\
\Hull(C_n).$\end{lemma}
\proof It is clear that $H$ contains the convex hull of $C_n.$ Since $H$
is closed, it contains the closure of $\cup\ \Hull(C_n).$ Moreover,
$\Hull(C_n)$ is an increasing sequence of convex sets and it follows that
$\cup\ \Hull(C_n)$ is convex.  Hence the closure is convex.  Thus $H$ is
contained in this closure.\endproof

\begin{proposition}
\label{hull_not_much_bigger}
Every point in the convex hull of ${\cal C}$ is within a distance $K_2$ of
${\cal C}.$
\end{proposition}
\proof The convex hull of a compact subset, $C,$ of ${\bf H}^3$ is equal to
the union of all geodesic tetrahedra having all four vertices in $C$ (we
allow degenerate tetrahedra where some of the vertices coincide).

Fix a point in $x$ in ${\cal C}$ and define $C_n$ to be the compact subset
of ${\cal C}$ of points within a distance $n$ of $x.$ Suppose that $z$ is a
point in $\Hull(C_n).$ By Lemma \ref{secondestimate}, $d(z,{\cal C})<K_2.$ Thus
$d(\Hull(C_n),{\cal C}) < K_2.$ Thus $d(\cup \ \Hull(C_n),{\cal C})\le K_2.$ Then
\ref{hullobservation} gives the result.  \endproof

\proof[Proof of \ref{Nproperty}]

By our initial careful choices of $N$, any translate of $\tilde{N}$ not
lying in ${\cal C}$ is very far from ${\cal C}$.  However, by
\ref{hull_not_much_bigger}, $\Hull({\cal C})$ is very close to ${\cal C}$,
whence the result.  \endproof

The authors thank A\thinspace W Reid  for carefully reading an early
version of this paper and the referee for several  useful comments.

Both authors are supported in part by the NSF.

\end{document}

%% file: gtoutput.tex
\def\ifplaintex{\expandafter\ifx\csname documentclass\endcsname\relax}

\ifplaintex 
\else
\fi

\expandafter\ifx\csname beginpicture\endcsname\relax
\expandafter\ifx\csname documentclass\endcsname\relax
\input pictex \else
\input prepictex \input pictex \input postpictex \fi\fi

\def\gt{{\mathsurround=0pt\it $\cal G\mskip-2mu$eometry \&\ 
$\cal T\!\!$opology}}        

\def\gtp{{\mathsurround=0pt\it $\cal G\mskip-2mu$eometry \&\ 
$\cal T\!\!$opology $\cal P\!$ublications}}  

\def\lognumber#1{\def\thelognumber{#1}}
\def\volumenumber#1{\def\thevolumenumber{#1}}
\def\papernumber#1{\def\thepapernumber{#1}}
\def\volumeyear#1{\def\thevolumeyear{#1}}

\def\pagenumbers#1#2{\def\startpage{#1}\def\finishpage{#2}}
\def\published#1{\def\publishdate{#1}}
\def\proposed#1{\def\theproposer{#1}}
\def\seconded#1{\def\theseconders{#1}}
\def\received#1{\def\receiveddate{#1}}
\def\revised#1{\def\reviseddate{#1}}
\def\accepted#1{\def\accepteddate{#1}}

\def\coverauthors#1{\def\thecoverauthors{#1}}
\def\asciiauthors#1{\def\theasciiauthors{#1}}

\long\def\asciiabstract#1{\long\def\theasciiabstract{#1}}

\let\\\par\let\thelognumber\relax
\let\thevolumenumber\relax\let\thepapernumber\relax
\let\thevolumeyear\relax\let\thesamplenumber\relax\let\startpage\relax
\let\finishpage\relax\let\publishdate\relax\let\receiveddate\relax
\let\reviseddate\relax\let\accepteddate\relax\let\theasciititle\relax
\let\theasciiauthors\relax
\let\theasciiabstract\relax
\let\theasciiemail\relax\let\theshortauthors\relax\let\theshorttitle\relax
\let\thecoverauthors\relax

\long\def\maketitlep{   

\count0=\startpage

\gt\hfill      
\beginpicture
\setcoordinatesystem units <0.33truein, 0.33truein> point at 2.2 0.9
\setplotsymbol ({$\cal G$})
\plotsymbolspacing=9truept
\circulararc 315 degrees from 0 1 center at 0 0
\setplotsymbol ({$\cal T$})
\circulararc 315 degrees from 1 -1 center at 1 0
\endpicture
\break
{\small\ifx\thesamplenumber\relax 
Volume \else Sample
\fi\thevolumenumber\ (\thevolumeyear)
\startpage--\finishpage\nl
Published: \publishdate}
\vglue 0.5truein plus 0.4fil minus 0.1truein

{\parskip=0pt\leftskip 0pt plus 1fil\def\\{\par\smallskip}{\ifplaintex\large
\else\Large\fi\bf\thetitle}\par\medskip}   

\vglue 0pt plus 0.1fil 

{\parskip=0pt\leftskip 0pt plus 1fil\def\\{\par}{\sc\theauthors}
\par\medskip}

\vglue 0pt plus 0.1fil 

{\small\parskip=0pt\let\newline\\
{\leftskip 0pt plus 1fil\def\\{\par}{\sl\theaddress}\par}
\expandafter\ifx\theemail\relax    
\relax\else\vglue 5pt plus 0.02fil minus 2pt\def\\{\stdspace{\rm 
and}\stdspace} 
\cl{Email:\stdspace\tt\theemail}\fi
\ifx\theurl\relax                  
\relax\else\vglue 5pt plus 0.02fil minus 2pt\def\\{\stdspace{\rm 
and}\stdspace}
\cl{URL:\stdspace\tt\theurl}\fi\par}

\vglue 7pt plus 0.3fil minus 3pt

{\bf Abstract}
\vglue 5pt plus 0.1fil minus 2pt

\theabstract

\vglue 7pt plus 0.3fil minus 3pt

{\bf AMS Classification numbers}\quad Primary:\quad \theprimaryclass

Secondary:\quad \thesecondaryclass

\vglue 5pt plus 0.3fil minus 2pt

{\bf Keywords}\quad \thekeywords

\vglue 10pt plus 0.5fil minus 5pt

{\small  Proposed: \theproposer\hfill Received: \receiveddate\nl
Seconded: \theseconders\hfill 
\ifx\reviseddate\relax                         
Accepted: \accepteddate                        
\else
Revised: \reviseddate                          
\fi}
\eject
}       

\let\maketitlepage\maketitlep
\let\maketitle\maketitlepage

\font\phead=cmsl9 scaled 950
\font=cmsl9 scaled 1050
\font\pnum=cmbx10 scaled 913
\font\lnum=cmbx10 
\font\pfoot=cmsl9 scaled 950
\font=cmsl9 scaled 1050
\ifplaintex
\headline{\vbox to 0pt{\vskip -4.5mm\line{\small\phead\ifnum
\count0=\startpage ISSN 1364-0380 (on line)
1465-3060 (printed) \hfill {\pnum\folio}\else\ifodd\count0\def\\{ }%
\ifx\theshorttitle\relax\thetitle\else\theshorttitle\fi\hfill{\pnum\folio}
\else\def\\{ and }{\pnum\folio}\hfill\ifx\theshortauthors\relax\theauthors
\else\theshortauthors\fi\fi\fi}\vss}}
\footline{\vbox to 0pt{\vglue 0mm\line{\small\pfoot\ifnum\count0=\startpage
\copyright\ \gtp\hfill\else
\gt, Volume \thevolumenumber\ (\thevolumeyear)\hfill\fi}\vss
}}
\else
\makeatletter
\def\@oddhead{{\small\lhead\ifnum\count0=\startpage ISSN 1364-0380 (on line)
1465-3060 (printed) \hfill {\lnum\number\count0}\else\ifodd\count0
\def\\{ }\ifx\theshorttitle\relax \thetitle \else\theshorttitle\fi\hfill
{\lnum\number\count0}\else\def\\{ and }{\lnum\number\count0}
\hfill\ifx\theshortauthors\relax 
\theauthors\else\theshortauthors\fi\fi\fi}}\def\@evenhead{\@oddhead}
\def\@oddfoot{\small\lfoot\ifnum\count0=\startpage\copyright\ \gtp\hfill\else
\gt, Volume \thevolumenumber\ (\thevolumeyear)\hfill\fi}
\def\@evenfoot{\@oddfoot}
\makeatother
\fi

\newwrite\gtoutfile
\long\gdef\makeheadfile{  
{\def\\{, }\def\s{ }
\immediate\openout\gtoutfile head.xxx
\immediate\write\gtoutfile{To: math@arxiv.org}
\immediate\write\gtoutfile{Subject: put or rep NNNNN:pppp}
\immediate\write\gtoutfile{--text follows this line--}
\immediate\write\gtoutfile{Proxy-for: \ifx\theasciiauthors\relax
\theauthors\else\theasciiauthors\fi\s<\ifx\theasciiemail\relax\theemail\else\theasciiemail\fi>}
\immediate\write\gtoutfile{\noexpand\\}
\immediate\write\gtoutfile{Authors: \ifx\theasciiauthors\relax
\theauthors\else\theasciiauthors\fi}
\immediate\write\gtoutfile{Title: \ifx\theasciititle\relax
\thetitle\else\theasciititle\fi}
\immediate\write\gtoutfile{Subj-class: GT or SG or MG etc}
\immediate\write\gtoutfile{MSC-class: \theprimaryclass\ifx\thesecondaryclass\relax\else, \thesecondaryclass\fi}
\immediate\write\gtoutfile{Journal-ref: Geom. Topol. \thevolumenumber
(\thevolumeyear) \startpage-\finishpage}
\immediate\write\gtoutfile{Comments: Published by Geometry and Topology at}
\immediate\write\gtoutfile{\s\s http://www.maths.warwick.ac.uk/gt/GTVol\thevolumenumber/paper\thepapernumber.abs.html}
\immediate\write\gtoutfile{\noexpand\\}
\immediate\write\gtoutfile{}
\ifx\theasciiabstract\relax
\immediate\write\gtoutfile{\theabstract}\else
\immediate\write\gtoutfile{\theasciiabstract}\fi
\immediate\write\gtoutfile{}
\immediate\write\gtoutfile{\noexpand\\}
\immediate\write\gtoutfile{}
\immediate\closeout\gtoutfile}}  

\def\maketitlepage{\maketitlep\makeheadfile}
\let\maketitle\maketitlepage

\def\ifplaintex{\expandafter\ifx\csname documentclass\endcsname\relax}

\ifplaintex 
\else
\fi

\expandafter\ifx\csname beginpicture\endcsname\relax
\expandafter\ifx\csname documentclass\endcsname\relax
\input pictex \else
\input prepictex \input pictex \input postpictex \fi\fi

\def\gt{{\mathsurround=0pt\it $\cal G\mskip-2mu$eometry \&\ 
$\cal T\!\!$opology}}        

\def\gtp{{\mathsurround=0pt\it $\cal G\mskip-2mu$eometry \&\ 
$\cal T\!\!$opology $\cal P\!$ublications}}  

\def\lognumber#1{\def\thelognumber{#1}}
\def\volumenumber#1{\def\thevolumenumber{#1}}
\def\papernumber#1{\def\thepapernumber{#1}}
\def\volumeyear#1{\def\thevolumeyear{#1}}

\def\pagenumbers#1#2{\def\startpage{#1}\def\finishpage{#2}}
\def\published#1{\def\publishdate{#1}}
\def\proposed#1{\def\theproposer{#1}}
\def\seconded#1{\def\theseconders{#1}}
\def\received#1{\def\receiveddate{#1}}
\def\revised#1{\def\reviseddate{#1}}
\def\accepted#1{\def\accepteddate{#1}}

\def\coverauthors#1{\def\thecoverauthors{#1}}
\def\asciiauthors#1{\def\theasciiauthors{#1}}

\long\def\asciiabstract#1{\long\def\theasciiabstract{#1}}

\let\\\par\let\thelognumber\relax
\let\thevolumenumber\relax\let\thepapernumber\relax
\let\thevolumeyear\relax\let\thesamplenumber\relax\let\startpage\relax
\let\finishpage\relax\let\publishdate\relax\let\receiveddate\relax
\let\reviseddate\relax\let\accepteddate\relax\let\theasciititle\relax
\let\theasciiauthors\relax
\let\theasciiabstract\relax
\let\theasciiemail\relax\let\theshortauthors\relax\let\theshorttitle\relax
\let\thecoverauthors\relax

\long\def\maketitlep{   

\count0=\startpage

\gt\hfill      
\beginpicture
\setcoordinatesystem units <0.33truein, 0.33truein> point at 2.2 0.9
\setplotsymbol ({$\cal G$})
\plotsymbolspacing=9truept
\circulararc 315 degrees from 0 1 center at 0 0
\setplotsymbol ({$\cal T$})
\circulararc 315 degrees from 1 -1 center at 1 0
\endpicture
\break
{\small\ifx\thesamplenumber\relax 
Volume \else Sample
\fi\thevolumenumber\ (\thevolumeyear)
\startpage--\finishpage\nl
Published: \publishdate}
\vglue 0.5truein plus 0.4fil minus 0.1truein

{\parskip=0pt\leftskip 0pt plus 1fil\def\\{\par\smallskip}{\ifplaintex\large
\else\Large\fi\bf\thetitle}\par\medskip}   

\vglue 0pt plus 0.1fil 

{\parskip=0pt\leftskip 0pt plus 1fil\def\\{\par}{\sc\theauthors}
\par\medskip}

\vglue 0pt plus 0.1fil 

{\small\parskip=0pt\let\newline\\
{\leftskip 0pt plus 1fil\def\\{\par}{\sl\theaddress}\par}
\expandafter\ifx\theemail\relax    
\relax\else\vglue 5pt plus 0.02fil minus 2pt\def\\{\stdspace{\rm 
and}\stdspace} 
\cl{Email:\stdspace\tt\theemail}\fi
\ifx\theurl\relax                  
\relax\else\vglue 5pt plus 0.02fil minus 2pt\def\\{\stdspace{\rm 
and}\stdspace}
\cl{URL:\stdspace\tt\theurl}\fi\par}

\vglue 7pt plus 0.3fil minus 3pt

{\bf Abstract}
\vglue 5pt plus 0.1fil minus 2pt

\theabstract

\vglue 7pt plus 0.3fil minus 3pt

{\bf AMS Classification numbers}\quad Primary:\quad \theprimaryclass

Secondary:\quad \thesecondaryclass

\vglue 5pt plus 0.3fil minus 2pt

{\bf Keywords}\quad \thekeywords

\vglue 10pt plus 0.5fil minus 5pt

{\small  Proposed: \theproposer\hfill Received: \receiveddate\nl
Seconded: \theseconders\hfill 
\ifx\reviseddate\relax                         
Accepted: \accepteddate                        
\else
Revised: \reviseddate                          
\fi}
\eject
}       

\let\maketitlepage\maketitlep
\let\maketitle\maketitlepage

\font\phead=cmsl9 scaled 950
\font=cmsl9 scaled 1050
\font\pnum=cmbx10 scaled 913
\font\lnum=cmbx10 
\font\pfoot=cmsl9 scaled 950
\font=cmsl9 scaled 1050
\ifplaintex
\headline{\vbox to 0pt{\vskip -4.5mm\line{\small\phead\ifnum
\count0=\startpage ISSN 1364-0380 (on line)
1465-3060 (printed) \hfill {\pnum\folio}\else\ifodd\count0\def\\{ }%
\ifx\theshorttitle\relax\thetitle\else\theshorttitle\fi\hfill{\pnum\folio}
\else\def\\{ and }{\pnum\folio}\hfill\ifx\theshortauthors\relax\theauthors
\else\theshortauthors\fi\fi\fi}\vss}}
\footline{\vbox to 0pt{\vglue 0mm\line{\small\pfoot\ifnum\count0=\startpage
\copyright\ \gtp\hfill\else
\gt, Volume \thevolumenumber\ (\thevolumeyear)\hfill\fi}\vss
}}
\else
\makeatletter
\def\@oddhead{{\small\lhead\ifnum\count0=\startpage ISSN 1364-0380 (on line)
1465-3060 (printed) \hfill {\lnum\number\count0}\else\ifodd\count0
\def\\{ }\ifx\theshorttitle\relax \thetitle \else\theshorttitle\fi\hfill
{\lnum\number\count0}\else\def\\{ and }{\lnum\number\count0}
\hfill\ifx\theshortauthors\relax 
\theauthors\else\theshortauthors\fi\fi\fi}}\def\@evenhead{\@oddhead}
\def\@oddfoot{\small\lfoot\ifnum\count0=\startpage\copyright\ \gtp\hfill\else
\gt, Volume \thevolumenumber\ (\thevolumeyear)\hfill\fi}
\def\@evenfoot{\@oddfoot}
\makeatother
\fi

\newwrite\gtoutfile
\long\gdef\makeheadfile{  
{\def\\{, }\def\s{ }
\immediate\openout\gtoutfile head.xxx
\immediate\write\gtoutfile{To: math@arxiv.org}
\immediate\write\gtoutfile{Subject: put or rep NNNNN:pppp}
\immediate\write\gtoutfile{--text follows this line--}
\immediate\write\gtoutfile{Proxy-for: \ifx\theasciiauthors\relax
\theauthors\else\theasciiauthors\fi\s<\ifx\theasciiemail\relax\theemail\else\theasciiemail\fi>}
\immediate\write\gtoutfile{\noexpand\\}
\immediate\write\gtoutfile{Authors: \ifx\theasciiauthors\relax
\theauthors\else\theasciiauthors\fi}
\immediate\write\gtoutfile{Title: \ifx\theasciititle\relax
\thetitle\else\theasciititle\fi}
\immediate\write\gtoutfile{Subj-class: GT or SG or MG etc}
\immediate\write\gtoutfile{MSC-class: \theprimaryclass\ifx\thesecondaryclass\relax\else, \thesecondaryclass\fi}
\immediate\write\gtoutfile{Journal-ref: Geom. Topol. \thevolumenumber
(\thevolumeyear) \startpage-\finishpage}
\immediate\write\gtoutfile{Comments: Published by Geometry and Topology at}
\immediate\write\gtoutfile{\s\s http://www.maths.warwick.ac.uk/gt/GTVol\thevolumenumber/paper\thepapernumber.abs.html}
\immediate\write\gtoutfile{\noexpand\\}
\immediate\write\gtoutfile{}
\ifx\theasciiabstract\relax
\immediate\write\gtoutfile{\theabstract}\else
\immediate\write\gtoutfile{\theasciiabstract}\fi
\immediate\write\gtoutfile{}
\immediate\write\gtoutfile{\noexpand\\}
\immediate\write\gtoutfile{}
\immediate\closeout\gtoutfile}}  

\def\maketitlepage{\maketitlep\makeheadfile}
\let\maketitle\maketitlepage

%% file: 2001-12.bbl
\begin{thebibliography}

\bibitem{B} {\bf F Bonahon}, {\em Bouts des vari\'{e}t\'{e}s hyperboliques de
dimension 3}, Ann.  of Math.  124 (1986) 71--158

\bibitem{BM} {\bf H Bass}, {\bf J Morgan}, {\em The Smith
Conjecture}, Academic Press (1984)

\bibitem{BH} {\bf S Bleiler}, {\bf C Hodgson}, {\em Spherical Space
Forms and Dehn Fillings}, Topology, 3 (1996) 809--833

\bibitem{CEG} {\bf D Canary}, {\bf D\,B\,A Epstein}, {\bf P Green},
{\em Notes on notes of Thurston}, from ``Analytic and Geometric
Aspects of Hyperbolic Space'', (D\,B\,A Epstein, editor) LMS Lecture
Notes, vol. 111 (1987)

\bibitem{CL} {\bf D Cooper}, {\bf D\,D Long}, {\em Virtually Haken
 surgery on knots},  J. Diff Geometry, 52 (1999) 173--187

\bibitem{CLR} {\bf D Cooper}, {\bf D\,D Long}, {\bf A\,W Reid}, {\em
Essential closed surfaces in bounded $3$--manifolds}, J.
Amer. Math. Soc.  10 (1997) 553--564

\bibitem{CG} {\bf M Culler}, {\bf C\,McA Gordon}, {\bf J Luecke},
{\bf P\,B Shalen},  {\em Dehn surgery on knots}, Ann.  of Math.
125 (1987) 237--300

\bibitem{dC} {\bf M do Carmo}, {\em Riemannian manifolds}, Birkhauser (1992)

\bibitem{F} {\bf S Fenley}, {\em Quasifuchsian Seifert Surfaces},
Math.  Zeit.  228 (1998) 221--227


\bibitem{He} {\bf J Hempel}, {\em 3--Manifolds}, Annals of Math Studies, 86,
Princeton University Press (1976) 


\bibitem{Ja}{\bf  W Jaco}, {\em Lectures on Three-Manifold Topology}, CBMS
series 43 (1977)

\bibitem{Li} {\bf T Li}, {\em Immersed essential surfaces in hyperbolic 
3--manifolds},  Comm. Anal. Geom. (to appear)



\bibitem{T}{\bf W\,P Thurston}, {\em The geometry and topology of
3--manifolds}, Princeton University Press

\end{thebibliography}
